\documentclass[11pt,a4paper]{amsart}
\usepackage{amsmath,amssymb}
\usepackage{bm}
\usepackage{graphicx}
\usepackage{ascmac}
 \usepackage{amsthm}
 
\theoremstyle{definition}
\newtheorem{thm}{Theorem}[section]
\newtheorem{dfn}[thm]{Definition}
\newtheorem{eg}[thm]{Example}
\newtheorem{cor}[thm]{Corollary}
\newtheorem{prop}[thm]{Proposition}
\newtheorem{lmm}[thm]{Lemma}
\newtheorem{rmk}[thm]{Remark}

\def\Supp{\mathop{\mathrm{Supp}}\nolimits}
\def\Spec{\mathop{\mathrm{Spec}}\nolimits}
\def\Proj{\mathop{\mathrm{Proj}}\nolimits}
\def\an{\mathop{\mathrm{an}}\nolimits}
\def\dd^c{\mathop{\mathrm{dd^c}}\nolimits}
\def\Res{\mathop{\mathrm{Res}}\nolimits}

\newcommand{\pr}{\mathbb{P}^1}
\newcommand{\pran}{\mathbb{P}^{1, \an}}

\newcommand{\cpx}{\mathbb{C}}
\newcommand{\real}{\mathbb{R}}
\newcommand{\aff}{\mathbb{A}^1}
\newcommand{\affan}{\mathbb{A}^{1,\an}}

\newcommand{\M}{\mathcal{M}}
\newcommand{\A}{\mathcal{A}}

\title{Activity measures of dynamical systems over non-archimedean fields}
\author{Reimi Irokawa}
\date{31/8/2019}
\address{Graduate School of Science, Tokyo Insitute of Technology}
\email{irokawa.r.aa@m.titech.ac.jp}
\subjclass[2010]{Primary: 37P50, Secondary: 14G22, 26E30, 37P30}
\keywords{non-archimedean dynamics; bifurcation; Berkovich geometry}
\begin{document}
\maketitle
\begin{abstract}
Toward the understanding of bifurcation phenomena of dynamics on the Berkovich projective line $\mathbb{P}^{1,an}$ over non-archimedean fields, we study the stability (or passivity) of critical points of families of polynomials parametrized by analytic curves. We construct the activity measure of a critical point of a family of rational functions, and study its properties. For a family of polynomials, we study more about the activity locus such as its relation to boundedness locus, i.e., the Mandelbrot set, and to the normality of the sequence of the forward orbit.
\end{abstract}
\section{Introduction}

Let $K$ be an algebraically closed field with complete, non-trivial and non-archimedean valuation. Let us consider an analytic family of rational functions $f : V\times_{K} \pran\to\pran$ of degree $d\geq2$ parametrized by a smooth strictly $K$-analytic curve $V$ and a marked point $c: V\to\pran$. The purpose of this paper is to develop the stability theory for such a pair $(f,c)$.

In complex dynamics, the stability, or the \textit{passivity} of the pair $(f,c)$ is well-studied in relation to the so-called J-stability, i.e., the stability of motion of Julia sets, and to the unlikely intersection problems in arithmetic dynamics. In non-archimedean dynamics, however, there are several difficulties even in giving an appropriate definition of the stability of $(f,c)$, which mainly comes from the lack of Montel’s theorem, although there are already several results by Silverman \cite{TS} and Lee \cite{Lee}. In this paper, we consider the stability from the following two approaches, and then study the relation between them, and moreover the relation between the stability in our sense and the J-stability.

The first approach is by means of normality. In complex dynamics, the passivity of the pair is originally defined in this way; namely, for a complex manifold $V$, a family $f:V\times\pr(\cpx)\to\pr(\cpx)$ of rational functions of degree $d$ and a marked point $c:V\to\pr(\cpx)$, the pair $(f,c)$ is \textit{passive} around $t_0\in V$ if the family $\{f_t^n(c(t))\}_n$ is normal around $t_0$, and \textit{active} if not. There is a fundamental result in complex dynamics by Ma\~ne-Sud-Sullivan \cite{MSS} and independently by Lyubich \cite{L83}, which states the following: for  $V$ and $f$ as above, with an additional assumption that every critical point moves analytically, $f$ is J-stable around $t_0$ if and only if all the critical points are passive around $t_0$. To show this, Montel’s theorem plays the essential role. In non-archimedean dynamics, however, the lack of Montel’s theorem is a big obstacle in the studies of stability theory. As a remedy, Favre-Kiwi-Trucco \cite{FKT} gave a weaker statement of Montel’s theorem by modifying the definition of normality, which we adopt to define the passivity of $(f,c)$; the pair $(f,c)$ is passive around $t_0\in V$ if $\{f_t^n(c(t)\}_n$ is normal around $t_0$ in the sense of \cite{FKT}. For details, see Definition \ref{FKTdef}. We treat this formulation in Section \ref{Normal}.

Our second approach is by means of the so-called \textit{activity measure}. The activity measure is introduced in complex dynamics by DeMarco \cite{DM1} to give a new condition for passivity of the pair $(f,c)$. Here, the activity measure, or the activity current for a higher-dimensional parameter space, is a locally-finite positive measure (resp. closed current) which can be written as the Laplacian of some continuous (pluri-)subharmonic function. DeMarco showed that the support of the activity current coincides with the activity locus of $(f,c)$. The importance of the activity measure comes not only from this fact but also from their nice potential-theoretic properties such as equidistribution as in \cite{DF}. In non-archimedean situation, on the other hand, we have good potentilal theory over smooth curves by Thuillier \cite{Thu}, whose discussion is almost parallel to that in complex case. By means of this, we will construct the non-archimedean activity measure $\mu_{(f,c)}$ in Section \ref{escaperate}; in case $f_t$ are polynomials, it can be written as the Laplacian of the function
\begin{equation}
h_{(f,c)}(t)=\lim_{n\to \infty}\frac{1}{d^n}\log\max\big(1,|f_t^n(c(t))|\big) \mbox{ on }V.
\end{equation}

Our activity measure for $(f,c)$ has the following properties:

\begin{thm}[Theorem \ref{3-main}, Proposition \ref{main2}]\label{equid}
{\it Let $K$, $V$, $f$ and $c$ be as above. We set $c_n(t)\colon=f_t^n(c(t))$. The sequence of positive measures $\{ \frac{1}{d^n}c_n^*\delta_{\zeta_{0,1}}\}_n$ converges weakly to the activity measure $\mu_{(f,c)}$, where $\delta_{\zeta_{0,1}}$ is the Dirac mass at the Gauss point $\zeta_{0,1}$ in $\pran$. Moreover, any exceptional point, i.e., a point in the set
\[
\mathcal{E}:=\bigg\{x\in\pran \bigg| \frac{1}{d^n}c_n^*(\delta_x-\delta_{\zeta_{0,1}})\not\rightarrow 0 \mbox{ weakly as }n\rightarrow\infty\bigg\}.
\]
is a type $1$ point.}
\end{thm}

Having these two approaches, we next discuss the interrelations between them, and compare them with the notion of J-stability. To this end, we first introduce the boundedness locus for a certain algebraic family of monic polynomials with $V=\affan$; when the characteristic of $K$ is $0$ or greater than $d$, we have a good parametrization of polynomials by critical points and the value at $0$ (for details, see Section \ref{poly}). In this case, we can consider the boundedness locus for the pair $(f,c)$ with $c$ critical:
\[
\M_{(f,c)}=\{t\in V|\{c_n(t)\} \mbox{ is bounded}\}.
\]
Assuming $\M_{(f,c)}$ is non-empty and bounded, we can show the following:

\begin{thm}[Theorem \ref{slicemain}]\label{Mandelbrot}
{\it Let $V$, $f$, $c$ and $\M_{(f,c)}$ be as above (see Section \ref{poly} for the precise definition). Then, the activity measure $\mu_{(f,c)}$ is the equilibrium measure of the set $\mathcal{M}_{(f,c)}$ with respect to $\infty$. In particular, the support of the activity measure coincides with the boundary of $\mathcal{M}_{(f,c)}$.}
\end{thm}

For the definition of the equilibrium measure, see Proposition \ref{defofequilm}. By using this property, we can compare the above two definitions of activity, and show a relation between passivity and J-stability by Silverman \cite{TS}.

\begin{prop}[Proposition \ref{main3}]\label{normalitythm}
\textit{Consider the same situation as in Theorem \ref{Mandelbrot} and assume that its Mandelbrot set $\M_{(f,c)}$ is non-empty and bounded. Then the passivity locus coincides with the normality locus.}
\end{prop}

\begin{cor}[Proposition \ref{main4}]\label{jstabilitythm}
{\it Let $K$, $V(=\affan)$, and $f$ be as in Proposition \ref{normalitythm}. Assume that the residue characteristic of $K$ is $0$ or greater than $d$. Then, every parameter $t_0$ where $f$ has an unstably indifferent periodic point lies in the activity locus of some critical point.}
\end{cor}

Silverman \cite{TS} showed that $f$ is J-stable around $t_0$ if and only if there exists a neighborhood $U$ of $t_0$ such that $f_t$ has no indifferent periodic point or no type $1$ repelling periodic point of multiplicity greater than $1$ for all $t$. The above corollary, therefore, indicates a partial relation between the J-stability and the stability of $(f,c)$.

\subsection{Organization of this paper}
In Section $2$, we collect basic results on potential theory over Berkovich curves. In Section $3$, we construct the activity measure and show an equidistribution result. In Section $4$, we consider some family of polynomials and discuss a certain relation between the activity measure and the boundedness locus. In Section $5$, we treat the passivity defined by the normality and consider the relation between the two notions of activity, and moreover, the relation between them and the J-stability. We give an example in Section $6$.

\subsection{Notations and conventions}
Throughout the paper, $K$ is an algebraically closed field with a complete, non-trivial and non-archimedean valuation. We denote $\underset{K}{\times}$ the fibred product and $\times$ the set-theoretic product. A smooth strictly $K$-analytic curve is a Berkovich analytic space which is smooth, good, paracompact, separated and purely of dimension $1$. 
%When we consider the projective line $\pran$, we always fix the infinity $\infty$ to consider the affine line $\pran=\affan\cup\{\infty\}$ in $\pran$, and the coordinate $x$ or sometimes $y$ on the affine line $\aff$. When we consider an affine line as parameter space, the coordinate is denoted by $t$.

\subsection*{Acknowledgements}

I would like to thank my advisor Fumiharu Kato for constant support and helpful conversations. Also I would like to thank Professor Mattias Jonsson and Professor Yutaka Ishii for giving me advices on this paper and my researches. This research has been partially supported by RIKEN Junior Research Associate Program and JSPS Grant-in-Aid for JSPS Fellows, 20J14309.

\section{Potential theory over smooth strictly $K$-analytic curves}\label{secpot}

\subsection{Potential functions and finite Radon measures}
For functions over a smooth and strictly $K$-analytic curve $V$, we can consider the Laplacian operator $\dd^c$. More precisely, the Laplacian operator is an $\mathbb{R}$-linear map from the space of functions of bounded differential variation to that of finite Radon measures. To study the weak convergence of some given sequence of finite Radon measures, it is useful to consider the sequences of the potential functions corresponding to the Radon measures.

\begin{prop}[\cite{Thu} Th\'eor\`eme 3.3.13.]\label{potmes}
\textit{ Assume that the variety $V$ is proper and irreducible. Then, for any Radon probability measure $\mu$ and any point $x\in V$ such that $x$ is not classical and does not belong to the support $\Supp\mu$ of $\mu$, there exists the unique function $u_{x,\mu}$ on $V$ such that
\[
u_{x,\mu}(x)=0, \mbox{ and}
\]\[
\dd^c u_{x,\mu}=\delta_x-\mu,
\]
where $\delta_x$ is the Dirac mass at the point $x$.}
\end{prop}

\begin{rmk}
%In the above theorem,``classical points" are $K$-rational ones on $V$.

For the definition of harmonic and subharmonic functions, see \cite{Thu}. These are completely analogous to the ones in complex potential theory.
\end{rmk}

\begin{eg}\label{egheight}
Let $z$ be an inhomogeneous coordinate of $\pr$. Let $h_{\infty}(z)$ be a continuous function defined by
\[
h_{\infty}(z)=\log\max(1,|z|).
\]
This function is non-constant only on the path from the point $\infty$ to the Gauss point $\zeta_{0,1}$. From this we have $\dd^c h_{\infty}(z)=\delta_{\infty}-\delta_{\zeta_{0,1}}$ and $h_{\infty}(\zeta_{0,1})=0$. The function $h_{\infty}$ is the potential function $-u_{\zeta_{0,1},\delta_{\infty}}$. We denote this function simply by $u_{\zeta_{0,1},\infty}$.
\end{eg}

Similarly, we denote the potential function $u_{x,\delta{y}}$ for a non-classical $x$ and any point $y$ simply by $u_{x,y}$. Later, we will regard them as two-variable function with respect to $x$ and $y$.

%\begin{eg}\label{pot}
%More generally, we can consider the function $u_{x,y}$ for any non-classical point $x$ and any point $y$. If there is a unique path $\gamma$ between these points, it is constant any path off from $\gamma$, and it is a non-constant affine function on $\gamma$.
%\end{eg}

\begin{eg}\label{Hsia}
%The following functions do not always satisfy the definition of potential functions, but still they are other important examples of functions on the analytification of the projective line $\pran$. 

In relation to the above functions, the following one is also important; for any $\zeta\in\pran$, the generalized \textit{ Hsia kernel} $\delta(x,y)_{\zeta}: \pran\times\pran\to\real$ with respect to $\zeta$ satisfies the following equation:
\[
\log\delta(x,y)_{\zeta}=u_{x,\zeta_{0,1}}(y)-u_{x,\zeta_{0,1}}(\zeta)- u_{y,\zeta_{0,1}}(\zeta).
\]
Since we know that $\dd^c u_{x,\zeta_{0,1}}(y)=\delta_{\zeta_{0,1}}-\delta_x$, we have
\[
\textstyle{\mathop{\mathrm{dd}}^c_x} \log\delta(x,y)_{\zeta}=\delta_{\zeta}-\delta_y,
\]
where $\mathop{\mathrm{dd}}^c_x$ is the Laplacian operator with respect to the variable $x$, which means we take the Laplacian as a function of $x$ while $y$ is considered to be constant. Since $\log\delta(x,\zeta)_{\zeta}=-u_{\zeta_{0,1},\zeta}(\zeta)$, the function $-\log\delta(x,y)_{\zeta}-u_{\zeta_{0,1},\zeta}(\zeta)$ is the potential function $u_{\zeta,y}(x)$ unless $\zeta$ is classical.

%We can represent the metric on $K=\affan(K)\subset\affan$ by means of the generalized Hsia kernel with respect to $\infty$; for all classical points $x$ and $y$, we have
%\[
%\delta(x,y)_{\infty}=|x-y|.
%\]
%Also, the spherical metric on the classical projective line $\pr(K)$ with a fixed homogeneous coordinate is a Hsia kernel with respect to $\zeta=\zeta_{0,1}$; for $x=[x_1:x_2]$ and $y=[y_1:y_2]$,
%\[
%\delta(x,y)_{\zeta_{0,1}}=\frac{|x_1y_2-x_2y_1|}{\max(|x_1|,|x_2|)\max(|y_1|,|y_2|)}.
%\]
%The Hsia kernel with respect to $\zeta_{0,1}$ is also called the spherical kernel.
\end{eg}

\subsection{Pull-back operators}
Let $W$ be another strictly $K$-analytic curve and $f:V\to W$ be a non-constant morphism. %In this section we define two operators on the space of positive Radon measures induced from $f$: the push-forward operator $f_*$ and the pull-back operator $f^*$. Since we only take pull-backs of measures, we just state the definition of pull back operators.  

\begin{dfn}
\textit{ In the above notation, let $x$ be a non-classical point of $W$. the pull-back $f^{*}\delta_{x}$ of the Dirac mass at $x$ by $f$ is defined to be the locally finite measure
\[
f^{*}\delta_{x}=\sum_{y\in f^{-1}\{x\}}m_{f}(y)\delta_y,
\]
where $m_{f}(y)$ is the multiplicity of $f$ at $y$.}
\end{dfn}

It is known in \cite{Thu} that any subharmonic function $h$ can be approximated by the decreasing sequence $\{h_n\}_n$ of smooth subharmonic functions, and the sequence of the Laplacians $\{\dd^c h_n\}$ converges to the Laplacian $\dd^c h$ of $h$. From this fact, together with Theorem \ref{potmes}, we can extend the pull-back operators to any finite positive Radon measures.

\begin{prop}[\cite{Thu}, Proposition 3.2.13]
\textit{ Let $V$ and $W$ be strictly $K$-analytic curves, $f\colon V\to W$ be a non-constant morphism, and $g$ be a smooth and subharmonic function on $W$. Then, we have
\[
\dd^c (g\circ f)=f^* \dd^cg.
\]}
\end{prop}
Note that, we can extend the same result to all the subharmonic function $g$, expressing it as a limit of smooth functions.

\begin{eg}[canonical measure]\label{canomeasure}
For a polynomial $\phi$ of degree $d$, we can consider the \textit{canonical height function} $h_{\phi}$ associated to $\phi$ given by
\[
h_{\phi}(x)=\lim_{n\to\infty}\frac{1}{d^n}\log\max(1,|\phi^n(x)|).
\]
This is a well defined, continuous and subharmonic function on $\pran$. Its Laplacian $\dd^c h_{\phi}$ is called the \textit{canonical measure}, whose support is called the Julia set of $\phi$.
\end{eg}
%For a complex rational function, we can also construct the canonical height function in similar way, but we need to take homogeneous lift and the argument is a little more cumbersome. The \textit{canonical measure} associated to $\phi$ is defined to be the Laplacian of the canonical height function. It is well-known that the support of this measure is exactly the Julia set of $\phi$. 

%In non-archimedean situation, we can also consider the canonical height function and the canonical measure, its Laplacian. All the above properties for the canonical heights and measures of complex rational functions, except for those about support set, are true for those of non-archimedean rational functions, too. The Julia set of $\phi$ is defined as the support of the canonical measure in non-archimedean setting. For details, see Chapter 10 of \cite{BR}. 

\subsection{Capacity theory for the projective line} \label{sec:capacity}

\begin{dfn}\textit{
On the Berkovich projective line $\pran$, for any probability measure $\mu$ and a point $\zeta$, the energy integral $I_{\zeta}(\mu)$ of $\mu$ with respect to $\zeta$ is defined to be
\[
I_{\zeta}(\mu)=\iint_{\pran\times \pran}-\log\delta(x,y)_{\zeta}(x)d\mu(x)d\mu(y).
\]
For a compact set $E\subset \pran\setminus \{\zeta\}$, let $\mathcal{P}(E)$ be the set of probability measures supported on $E$. The Robin constant $V_{\zeta}(E)$ of $E$ with respect to the point $\zeta$ is defined to be
\[
V_{\zeta}(E)=\inf_{\mu\in\mathcal{P}(E)}I_{\zeta}(\mu).
\]
The logarithmic capacity of $E$ with respect to the point $\zeta$ is
\[
C(E;\zeta)=e^{-V_{\zeta}(E)}.
\]}
\end{dfn}

\begin{prop}[\cite{BR}, Proposition 6.6 \& Proposition 7.21] \label{defofequilm}
\textit{ Assume that $C(E;\zeta)>0$. There exists an unique measure $\mu_E\in\mathcal{P}(E)$ such that
\[
V_{\zeta}(E)=I_{\zeta}(\mu_E).
\]}
\end{prop} 

We call $\mu_E$ the \textit{equilibrium measure} of $E$ with respect to the point $\zeta$. The potential function $u_{\zeta, \mu_E}$ is also denoted by $u_{\zeta,E}$. It might depend on the point $\zeta$, but actually, we have the following fact:

\begin{prop}
\textit{ If the compact set $E$ is of capacity $0$ with respect to some point $\zeta$, i.e.\ $C(E;\zeta)=0$, then it is of capacity $0$ with respect to any point $\zeta'$.}
\end{prop}

Hence, whenever we only care as to whether the capacity is zero, it is enough to discuss the capacity for some fixed point. We mainly consider the capacity with respect to $\infty$ in the later argument.

The next proposition is used later:

\begin{prop}[\cite{BR}, Proposition 6.8 \& Corollary 7.39]\label{contipot}
\textit{ In the above notation, assume that $E\subset V\setminus \{\zeta\}$ has positive capacity. Let $U_{\zeta}$ be a connected component of $V\setminus E$ containing $\zeta$. Then, the support of the equilibrium measure $\mu_E$ is contained in the boundary $\partial U_{\zeta}$ of $U_{\zeta}$. Moreover, if $u_{\zeta, E}$ is continuous, the support of the measure $\mu_E$ coincides with $\partial U_{\zeta}$.}
\end{prop}

\subsection{Arakelov-Green functions} \label{subsec:arakelovgreen}
%In this section, we only consider functions on the projective line $\pran$. 
On $\pran$, we consider the following property to discuss Arakelov-Green functions:
\begin{dfn}
\textit{ A probability measure $\mu$ is said to be of continuous potentials if for some non-classical point $\zeta$ the potential function $u_{\zeta,\mu}$ is continuous.}
\end{dfn}

Note that it can be easily seen that if $\mu$ is of continuous potentials, then the potential function $u_{\zeta,\mu}$ is continuous for any non-classical point $\zeta$.

%\begin{prop}
%\textit{For a probability measure $\mu$, the following conditions are equivalent:
%\begin{enumerate}
%\item for all non-classical point $\zeta$, the potential function $u_{\zeta,\mu}$ is continuous;
%\item for some non-classical point $\zeta$, the potential function $u_{\zeta,\mu}$ is continuous.
%\end{enumerate}}
%\end{prop}
%This can be easily seen by the fact the potential function $u_{\zeta,\xi}$ is continuous for any non-classical points $\zeta$ and $\xi$.

\begin{dfn}
\textit{For a probability measure $\mu$ with continuous potentials, an\textit{ Arakelov-Green function} $g_{\mu}\colon\pran\times\pran\to\real$ is defined to be
\begin{equation}\label{AG*}
g_{\mu}(x,y)=\int_{\pran} -\log\delta(x,y)_{\zeta}d\mu(\zeta)+ C, 
\end{equation}
where $C$ is some constant. }
\end{dfn}

Of course we can define $g_{\mu}(x,y)$ for any probability measure $\mu$ not necessarily with continuous potentials. The following properties, however, are important and valid only for ones with continuous potentials:
\begin{prop}[\cite{BR}, Proposition 8.66]
\textit{ For a probability measure $\mu$ with continuous potentials, the Arakelov-Green function $g_{\mu}(x,y)$ is lower-semicontinuous as a function of $2$-variables, continuous as a single-variable function of each variable, and symmetric with respect to $x$ and $y$. It is of bounded differential variation, and
\[
\dd^c_{x} g_{\mu}(x,y)=\delta_{y}-\mu.
\]}
\end{prop}
\begin{prop}[Energy-minimizing principle, \cite{BR}, Theorem 8.71]\label{EMP}
\textit{For a probability measure $\mu$ with continuous potentials, consider the energy integral
\[
I_{\mu}(\rho)=\iint_{V\times V} g_{\mu}(x,y)d\rho(x)d\rho(y),
\]
for any probability measure $\rho$. Then, the integral takes its minimum if and only if $\rho=\mu$.}
\end{prop}

\subsection{Normality of functions over Berkovich curves} \label{subsection:normality}
 
 In complex dynamics, the following Montel's theorem plays a fundamental role; for a complex manifold $V$, any family of analytic functions from $V$ to $\pr(\cpx)$ which avoid $0$, $1$ and $\infty$ is normal. Unfortunately, it is known that this theorem is not true in non-archimedean setting.

In \cite{FKT}, however, they give an alternative definition of normality so that a similar result for the Montel's theorem holds.

\begin{dfn}\label{FKTdef}
\textit{Let $X$ be an open subset of $\pran$. A family $\mathcal{F}$ of analytic functions on $X$ with values in $\pran$ is normal if for any sequence $f_n \in\mathcal{F}$ and any point $x\in X$, there exists a neighborhood $V$ of $x$, and a subsequence $\{f_{n_j}\}_j$ which converges pointwise on $V$ to a continuous function.}
\end{dfn}

%The following theorem is so-called the non-archimedean Montel's theorem:

%\begin{prop}[\cite{FKT}, Theorem 5.4]\label{FKTmain}
%\textit{Any family of meromorphic functions on an open subset $X$ of $\pran$ such that, for all $x\in X$, local unseparable degrees at $x$ are bounded, and avoids three points in $\pran$ is both normal, and equicontinuous at any classical point.}
%\end{prop}

We will use the next theorem, one of the consequences of the non-archimedean Montel's theorem;

\begin{prop}[\cite{FKT}, Theorem 2.1]\label{FKTthm}\textit{
    Let $U$ be a basic subset, i.e., the boundary of $U$ consists of only type $2$ or $3$ points, and $Y$ be an affinoid domain of $\pran$. Then any sequence of analytic functions $\{f_n:U\to Y\}$ admits a subsequence that converges pointwise on $U$ to a continuous function.
}\end{prop}

%% DYNAMICAL GREEN FUNCTIONS %%
\section{The construction of the activity measure}\label{escaperate}

Let $V$ be a smooth strictly $K$-analytic curve, $f: V\underset{K}{\times} \pran\rightarrow\pran$ be an analytic family of rational functions of degree $d$, and $c: V\rightarrow \pr$ be a marked point. %For any $t\in V$, $f_t :\{t\}\underset{K}{\times} \pran\to\{t\}\underset{K}{\times}\pran$ defines a rational functions over the analytification of the projective line over some extension of the valuation field $K$. We assume the degree of $f_t$ is $d$ for any $t$, which is the meaning of the analytic family of rational functions of degree $d$. 
The goal of this section is as follows:

\begin{thm} \label{3-main} \textit{
In the above settings, there exists a locally finite positive Radon measure $\mu_{(f,c)}$, locally written as the Laplacian of some continuous subharmonic function, which satisfies the following property: For the family of analytic functions $\{c_n(t):=f_t^n(c(t))\colon V\to\pran\}_n$, we have
\begin{align*}
    \frac{1}{d^n}c_n^*\delta_{0,1}\to \mu_{(f,c)} \text{ weakly as }n\to\infty.
\end{align*}
Moreover, if $f_t$ is a family of polynomials, the measure $\mu_{(f,c)}$ can be globally written as the Laplacian of a continuous and subharmonic function on $V$.
}
\end{thm}

\begin{dfn}\textit{
The measure $\mu_{(f,c)}$ defined above is called the activity measure of the pair $(f,c)$.
}
\end{dfn}

For a family of polynomials $f_t$ and its critical point $c\neq\infty$, Theorem \ref{3-main} is easy; define the function $h_{(f,c)}$ by
\begin{align*}
    h_{(f,c)}(t)=\lim_{n\to\infty}\frac{-1}{d^n}\log\max(|c_n(t)|,1).
\end{align*}
Then, clearly the limit $h_{(f,c)}$ exists, and it is naturally extended to a continuous subharmonic function. Since
\begin{align*}
    \mathop{\dd^c}(-\log\max(|c_n(t)|,1))&=c_n^*(\delta_{0,1}-\delta_{\infty}) \\
    &=c_n^*\delta_{0,1},
\end{align*}
we have the theorem for the families of polynomials.

For general families of rational functions, we take a homogeneous lift of $f_t$ and define the potential function locally, and then patch their Laplacians together to get the required measure defined on whole $V$.

For any point $t\in V$, take a sufficiently small affinoid neighbourhood $\M(\A)$ of $t$ in $V$ such that
\[
f_t(z)=\frac{q_t(z)}{p_t(z)} \mbox{ for some }p_t,\ q_t\in\A[z].
\]
Let $P_t(X,Y)=P_t^{(1)}(X,Y)$ and $Q_t=Q_t^{(1)}(X,Y)$ be a homogeneous lifts of $p_t$ and $q_t$ respectively with $z=Y/X$, with non-zero resultant $\Res(P_t,Q_t)\neq0$ for any $t\in\M(\A)$. Define $P_t^{(n)}=P_t^{(n-1)}(P_t,Q_t)$ and $P^{(1)}_t=P_t$, and $Q_t^{(n)}=Q_t^{(n-1)}(P_t,Q_t)$ and $Q^{(1)}_t=Q_t$ respectively. Note that we have $\Res(P_t^{(n)},Q_t^{(n)})\neq0$ for any $t$.

Next, let us consider a lift $C(t)\colon V\to \mathbb{A}^{2,an}\setminus\{0\}$ of $c(t)\colon V\to\pran$ with no pole as follows; for a sufficiently small neighborhood $\A$, we may assume $c(t)\in \affan$ for any $t$ by changing the coordinate if necessary. The lift is defined by $C(t)=(1,c(t))$.

%Set $\mathbb{A}_{0}=\pran\setminus\{\infty\}$, $\mathbb{A}_{\infty}=\pran\setminus\{0\}$, $V_0=c^{-1}(\mathbb{A}_{0})$ and $V_{\infty}=c^{-1}(\mathbb{A}_{\infty})$, where $c^{-1}$ is the pullback of $c$. The morphism $c$ can be regarded as a section on $V_0$ and $V_{\infty}$ respectively via $\mathbb{A}_{0}\simeq \affan$ and $\mathbb{A}_{\infty}\simeq \affan$. Therefore, the homogenous lift $C$ of $c$ is defined as
%\[
%C(t)=\begin{cases}
%(1, c(t)) \mbox{ if }t\in V_0, \\
%(c(t), 1) \mbox{ if }t\in V_{\infty}.
%\end{cases}
%\]
Now we construct locally the activity measure. Take any $t\in V$, an affinoid neighbourhood $\M(\A)$ of $t$, and an open neighborhood $U$ of $t$ in $\M(\A)$. Set
\begin{align*}
    h^{(n)}(t)=\frac{1}{n}\log\left\|P_t^{(n)}(C(t)),Q_t^{(n)}(C(t))\right\|,
\end{align*}
where $\|\cdot,\cdot\|$ is a generalized Hsia kernel with respect to $\zeta_{0,1}$; note that it is an upper-semicontinuous function of two variables satisfying $\|x,y\|=\max(|x|,|y|)$ for any $x,y\in K$.
%\[
%h^{(n)}(t)=\frac{1}{n}\log\max\left(\left|P_t^{(n)}(C(t))\right|,\left|Q_t^{(n)}(C(t))\right|\right) \mbox{ on }U.
%\]

\begin{prop}\textit{
In the above notation, for any $t\in U(K)$, we have positive functions $C_1(t)$ and $C_2(t)$ such that
\[
C_1(t)\|x,y\|^{d}\leq \|P_t(x,y),Q_t(x,y)\|\leq C_2(t)\|x,y\|^{d},
\]
for all $t\in U$ of type $1$ and all $(x,y)\in\mathbb{A}^2(K)$.
}\end{prop}
\begin{proof}
Write $P_t(z)=\sum_i a_i(t) z^i$ and $Q_t(z)=\sum_i b_i(t)z^i$. Then we have $C_2(t)=\max_{i,j}(|a_i(t)|,|b_j(t)|)$. For $C_1(t)$ we need a little more complicated argument. First of all, since $\Res(P_t,Q_t)\neq 0$ there exist polynomials $G_{1,t}$, $G_{2,t}$, $H_{1,t}$, and $H_{2,t}\in \mathbb{Z}[\{a_i(t)\},\{b_i(t)\},X,Y]$ homogeneous in $X$ and $Y$  of degree $d-1$ such that
\begin{align*}
P_tG_{1,t}+Q_tG_{2,t}&=\Res(P_t,Q_t)X^{2d-1},\mbox{ and} \\
P_tH_{1,t}+Q_tH_{2,t}&=\Res(P_t,Q_t)Y^{2d-1}.
\end{align*}
Taking the maximum of the absolute values of both sides of these equations, we have
\begin{align*}
    \max(|P_tG_{1,t}|,|Q_tG_{2,t}|,|PH_{1,t}|,|QH_{2,t}|)\geq|\Res(P_t,Q_t)|\max(|X|,|Y|).
\end{align*}

Since $G_{1,t}$, $G_{2,t}$, $H_{1,t}$, and $H_{2,t}\in \mathbb{Z}[\{a_i(t)\},\{b_i(t)\},X,Y]$ are homogeneous of degree $d-1$, we have
\begin{align*}
    |G_{1,t}(X,Y)|\leq \max(|a_i(t)|,|b_i(t)|)\cdot \max(|X|,|Y|)^{d-1},
\end{align*}
and the similar inequality holds for $G_{2,t}$, $H_{1,t}$ and $H_{2,t}$. Therefore,
\begin{align*}
    \max(|P_t|,|Q_t|)\max(|a_i(t)|,|b_i(t)|)\max(|X|,|Y|)^{d-1}\geq |\Res(P_t,Q_t)|\max(|X|,|Y|)^{2d-1}.
\end{align*}

For each $t$, taking $B(t)=\max(|a_i(t)|,|b_i(t)|)>0$, we have
\begin{align*}
    B(t)\|P_t(X,Y),Q_t(X,Y)\|\geq|\Res(P_t,Q_t)|\cdot\|X,Y\|.
\end{align*}
Therefore, by setting $C_1(t)=|\Res(P_t,Q_t)|/B(t)$ we have
\[
C_1(t)\|x,y\|^d\leq \|P_t(x,y),Q_t(x,y)\|.
\]
\end{proof}
\begin{prop}
\textit{In the above notation, the sequence $\{h^{(n)}\}$ converges uniformly to a continuous and subharmonic function (denoted by $h_{(f,c)}$).}
\end{prop}

\begin{proof}
%In particular,
Taking $(x,y)=C(t)$ we have
\begin{equation}\label{res}
C_1(t)\|C(t)\|^d\leq \|P_t(C(t)),Q_t(C(t))\|\leq C_1(t)\|C(t)\|^d
\end{equation}
for all $t\in U$ of type $1$. Since $C_1$ and $C_2$ are determined by taking maximum, multiplication, addition, and division of coefficients of $P_t$ and $Q_t$. the functions $\log C_1(t)$ and $\log C_2(t)$ can be extended to whole $U$ continuously. Then, since the set of type $1$ points is dense in $U$, we have the formula (\ref{res}) for all $t\in U$. Therefore, for each compact subset $E$ of $U$, letting $C=\log\max_t(C_1(t), C_2(t))$, we have
\[
\log\|C(t)\|\leq \frac{1}{d}\log\|P_t(C(t)),Q_t(C(t))\|+C,
\]
which shows the locally uniform convergence of $\{h^{(n)}\}$ on $E$.

Continuity and subharmonicity follow from those of $\{h^{(n)}\}$
\end{proof}
%Remark that when $f$ is a family of polynomials, we have
%\[
%h_{(f,c)}(t)=\lim_{n\to\infty}\frac{1}{d^n}\log\max(|f^n_t(c(t))|,1).
%\]
%In this case we don't have to take an open set $U$ of $V$ and $h_{(f,c)}$ is defined on whole $V$. 

Define the measure $\mu_{(f,c)}^U=-\dd^c h_{(f,c)}$ on $U$. This is a positive finite Radon measure on $U$.

\begin{prop}\textit{
In the above notations, define $c_n(t)=f_t^n(c(t))$. Then, we have
\[
\frac{1}{d^n}(c_n|_{U})^{*}\delta_{\zeta_{0,1}^U}\to\mu_{(f,c)} \mbox{ as }n\to\infty\text{ on }U,
\]
where the limit is weak convergence.
}\end{prop}

\begin{proof}

we have
\begin{eqnarray}\label{hgt}
h^{(n)}(t)&=&\frac{1}{d^n}\log\max(|P_t^{(n)}(C(t))|,|Q_t^{(n)}(C(t))|)\nonumber\\
&=&\frac{1}{d^n}\bigg(\log\max(|f^n_t(c_t)|,1)+\log\max(|P_t^{(n)}(C(t))|,1)+\log|c_1(t)|\bigg),
\end{eqnarray}
where $C(t)=(c_1(t),c_2(t))$ for all $t\in V\setminus\{t$ $\colon f_t(c(t))=\infty$ or $c(t)=\infty\}$. We can extend the equality if we take the limit for each point we exclude. As we see in Example \ref{egheight} we have
\[
\dd^c \log\left|P_t^{(n)}(C(t))\right|=c_n^*(\dd^c \delta_{\infty}-\delta_{\zeta_{0,1}}) \mbox{ on }U.
\]
For the third term of (\ref{hgt}), we have
\[
\dd^c\log|c_1(t)|=\sum_{x\in V(K)} a(x)\delta_{x},
\]
where div $c_1=\sum a(x)[x]$. We have the similar formula for the second one. By the hypothesis on the lifts of $f_t$ and $c_t$, we know that
\begin{itemize}
\item $c_2$ has no pole,
\item if $c_2$ has a zero, then $|P_t(c(t))|$ has a pole of the same degree
\item if $|P_t(C(t))|$ has a zero, then $|f_t(c(t))|$ has a pole of the same degree,
\item $|P_t(C(t))|$ and $|Q_t(C(t))|$ has no pole.
\end{itemize}
Therefore, we have
\[
\dd^c\frac{1}{d^n}\bigg(\log\max(|f_t(c_t)|,1)+\log\max(|c(t)|,1)+\log|c_1(t)|\bigg)=-\frac{1}{d^n}c_n^*\delta_{\zeta_{0,1}},
\]
where $c_n(t)=f_t^n(c(t))$. As $h^{(n)}\to h_{(f,c)}$ on $U$ locally uniformy as $n\to\infty$, we have 
\begin{equation}\label{main1}
\lim_{n\to\infty}\frac{1}{d^n}c_n^*\delta_{\zeta_{0,1}}=\mu_{(f,c)}^{U}\mbox{ weakly},
\end{equation}
which is a positive finite Radon measure.
\end{proof}

\begin{prop}\textit{
The above measure $\mu_{(f,c)}^U$ is independent of the choice of the open set $U$ and the lifts $P_t$, $Q_t$ and $C(t)$, i.e., if $U'$ is another open set and $P'_t$, $Q'_t$ and $C'(t)$ is a lift of $(f,c)$ on $U'$ satisfying the same condition as $P_t$, $Q_t$ and $C(t)$, then
\[
\mu_{(f,c)}^U|_{U\cap U'}=\mu_{(f,c)}^{U'}|_{U\cap U'}.
\]
}
\end{prop}
\begin{proof}
%Let $C'(t)=(C'_1(t),C'_2(t))\colon V\to \mathbb{A}^{2,an}$ be another lift of $c(t)$ which has no pole on $V$. Then,
Since both of $C(t)$ and $C'(t)$ have no zero or pole on $U$ and $U'$ respectively, there exists an invertible analytic function $\varphi$ on $U\cap U'$ such that
\[
C'(t)=\varphi(t)\cdot C(t) \text{ on }U\cap U'.
\]

For $(P_t, Q_t)$ and $(P'_t,Q'_t)$, since both of them represents the same rational function, we have
\begin{align*}
    P'_t(X,Y)&=a(t)P_t(X,Y)\text{, and} \\
    Q'_t(X,Y)&=a(t)Q_t(X,Y) \text{ on } U\cap U',
\end{align*}
where $a(t)$ is an invertible analytic function on $U\cap U'$. Note that since $\varphi$ and $a(t)$ are invertible, $\log|\varphi|$ and $\log|a|$ are harmonic. Therefore, we have
\begin{gather*}
\log\max\left(\left|P'^{(n)}_t(C'(t))\right|,\left|Q'^{(n)}_t(C'(t))\right|\right)\\ =\log\max\left(\left|P^{(n)}_t(C(t))\right|,\left|Q^{(n)}_t(C(t))\right|\right)+d^n\log|\varphi(t)|+(d^{n-1}+1)\log|a(t)|.
\end{gather*}
for any $n=1,2,3,\ldots$ on $U\cap U'$. Dividing them by $d^n$ and taking a limit, we have
\[
h_{(f,c)}(t)=h'_{(f,c)}(t)+\log|\varphi(t)|+\log|a(t)| \text{ on }U\cap U',
\]
where $h'_{(f,c)}(t)$ is defined in the same way as $h_{(f,c)}$ replacing 
$P_t$, $Q_t$, and $C$ by $P'_t$, $Q'_t$ and $C'$ respectively.
%respectively. 
 As a consequence, we have $\dd^c h_{(f,c)}=\dd^c h'_{(f,c)}$ on $U\cap U'$.
\end{proof}

The above proposition shows the compatibility of the system of locally defined measures $\{\mu_{(f,c)}^U\}$. Hence we can glue them together to get the global measure $\mu_{(f,c)}$ defined on the whole $V$.
%begin{dfn}\textit{
%The measure $\mu_{(f,c)}$ constructed above is called the activity measure of the pair $(f,c)$.}\end{dfn}

Since the formula (\ref{main1}) is valid locally, this holds globally on $V$. %We restate this result as a theorem:
%\begin{thm}\label{main1th}\textit{
%Let $V$ be an smooth strictly $K$-analytic curve, $f\colon V\times_K \pran\to\pran$ be an analytic family of rational functions of degree $d$, and $c$ be a marked point. Set $c_n(t)=f_t^n(c(t))$. Then, we have
%\[
%\frac{1}{d^n}c_n^{*}\delta_{\zeta_{0,1}}\to\mu_{(f,c)} \mbox{ as }n\to\infty,
%\]
%where the limit is weak convergence.
%}\end{thm}
This completes the proof of Theorem \ref{3-main}.

Note that the activity measure is defined independently of the choice of the Gauss point. More generally, we have the following:

\begin{prop}\label{main2}\textit{
For any non-classical point $\zeta\in\pran$, we have
\[
\frac{1}{d^n}c_n^{*}\delta_{\zeta}\to\mu_{(f,c)} \mbox{ as }n\to\infty,
\]
where the limit is weak convergence.
}\end{prop}
\begin{proof}
For any point $\zeta$, define the function $u_{\zeta}$ to be
\[
u_{\zeta}(x)=\log\delta(x,\zeta)_{\zeta_{0,1}}
\]
defined in Example \ref{Hsia}. This is a bounded and continuous function on $\pran$, whose Laplacian is $\dd^cu_{\zeta}=\delta_{\zeta_{0,1}}-\delta_{\zeta}$. For each $n$, consider $u_{\zeta}(c_n(t)):V\to\real$. This is also a bounded continuous function satisfying $\dd^c u_{\zeta}(c_n(t))=c_n^{*}(\delta_{\zeta_{0,1}}-\delta_{\zeta})$. Since
\[
\lim_{n\to\infty}\frac{1}{d^n}u_{\zeta}(c_n(t))=0
\]
uniformly on $V$, we have $\frac{1}{d^n}c_n^*(\delta_{\zeta_{0,1}}-\delta_{\zeta})\to 0$ weakly as $n\to\infty$.
\end{proof}

\begin{rmk}
The exceptional set 
\[
\mathcal{E}:=\big\{x\in\pran\ |\ \frac{1}{d}(\delta_x-\delta_{\zeta_{0,1}})\not\to 0\mbox{ as }n\to\infty\big\}
\]
is actually of capacity $0$. We show this fact in a forthcoming paper with Y.\ Okuyama.

%a subset of the set of type $1$ point $\pr(K)$ by the above proposition, but  is expected to be polar. So far nothing more is known. The difficulty of the evaluation of $\mathcal{E}$ comes from the non-existence of Lebesgue measure on $V$.
\end{rmk}

\section{The activity measure vs. the boundedness locus}\label{poly}
In this section, we will consider a relation between the activity measure and the boundedness locus for a certain family of polynomials.
\subsection{The boundedness locus}\label{poly-1}
In this section, we consider the following situation. We assume that the characteristic of $K$ is $0$ or greater than $d$. For $\textbf{c}=(c_1,\ldots,c_{d-2})\in K^{d-2}$ and $a\in K$, define a polynomial $f_{\textbf{c},a}$ of degree $d$ as follows:
\[
f_{\textbf{c},a}(z)=\frac{1}{d}z^d+\sum_{k=1}^{d-2}\frac{\sigma_{k}(c_1,\ldots,c_{d-2})}{d-k}z^{d-k}+a^d,
\]
where $\sigma_{k}$ is the $k$-th elementary symmetric polynomial. Then, the critical points of $f_{\textbf{c},a}$ are $c_0=0,c_1,\ldots,c_{d-2}$ and $f_{\textbf{c},a}(0)=a^d$. Regarding the set $\mathbb{A}^{d-2}_K\times_K\mathbb{A}_K$ of points $(\textbf{c},a)$ as a parameter space of polynomials, we consider its compactification $\mathbb{P}:=\mathbb{P}^{d-1}_K$ and its boundary $\mathbb{P}_{\infty}:=\mathbb{P}\setminus(\mathbb{A}^{d-2}_K\times_K\mathbb{A}_K)$. The family we consider in this section is a $1$-parameter subfamily of this $\mathbb{A}^{d-2}_K\times_K\mathbb{A}_K$ given as follows; take any line $\overline{V}$ of $\mathbb{P}$ such that $\overline{V}\cap\mathbb{P}_{\infty}$ consists of a single point denoted by $\infty$. Let $V$ denote the subset of the parameter space $\mathbb{A}^{d-2}\times_K\mathbb{A}$ defined by $V=\overline{V}\setminus\{\infty\}$, which we identify with $\aff$ by an affine coordinate $t$. By abuse of notation, we denote by $\overline{V}$ and $V$ their analytification respectively. Let $f_t=f_{\textbf{c}(t),a(t)}$ denote the analytification of the restriction to $V$ of the map $f_{\textbf{c},a}$. For a marked point, take one of $c_i$'s and denote by $c$ its analytifictaion. 
%We denote the parameter space by $\overline{V}$ and define $V:=\pr\setminus\{\infty\}$. Let $f_t$ denote the analytification of the subfamily $f_t=f_{(\textbf{c}(t),a(t)}$. For a marked point, take one of $c_i$'s or $0$ and denote $c$ by its analytification. 
We consider the relation between the activity and the boundedness of the orbit of $c$ by $f$ i.e.\ $\{c_n(t)\}_n$.

%Then we can consider the compactification $\mathbb{P}:=\mathbb{P}^{d-1}$ of $\mathbb{A}^{d-1}$ by $(\textbf{c},a)\mapsto[c_1:\cdots:c_{d-2}:a:1]$. We set $\mathbb{P}_{\infty}=\mathbb{P}\setminus\mathbb{A}^{d-1}(\simeq \mathbb{P}^{d-2})$. Consider any slice by $\pr$ of the compactification $\mathbb{P}$ of the moduli space such that a subspace $\pr\cap\mathbb{P}_{\infty}$ consists of a single point denoted by $\infty$. In this setting, we can consider the dynamical Green function $h_{(f,c)}$ where the parameter space $V$ is the analytification of $\pr\setminus\{\infty\}$ and $c$ is one of the critical points of $f_{c,a}$. 
We consider the boundedness locus of a given critical point, which is a higher-degree analogue of Mandelbrot set: for any pair $(f,c)$ considered above, define
\[
\mathcal{M}_{f,c}=\left\{ t\in V\middle| \{f_t^n(c(t))\}_{n=0}^{\infty}\text{ is bounded}\right\}.
\]

\begin{rmk}
 Consider the family of quadratic polynomials $f_t(z)=z^2+t$ ($t\in\affan$). In this case, the point $0$ is the unique critical point for any $t$ and we have the boundedness locus, i.e.\ the ``Mandelbrot set"
\[
\mathcal{M}=\{ t\in \affan| \{f_t^n(0)\}_{n=0}^{\infty}\text{ is bounded}\}.
\]
This set is, in non-arhimedean case, actually, just the closed unit disc (for details, see Section \ref{quad}). The above $\mathcal{M}_{f,c}$ is, hence, despite with slightly different parametrization, can be seen as some sort of higher degree generalization of the Mandelbrot set. 
\end{rmk}

From now on, we assume that the set $\mathcal{M}_{f,c}$ is a non-empty bounded set of the parameter space $V$. This condition is not extremely strong, as we discuss in Section \ref{bddcond}. In this case we have the following:
\begin{prop}\textit{
If $\mathcal{M}_{f,c}$ is non-empty, the activity measure $\mu_{(f,c)}$ is a strictly positive finite Radon measure.
}\end{prop}

\begin{proof}
The function $h_{(f,c)}$ is defined as a non-constant function on $V$ since $h_{(f,c)}(t)=0$ for any $t$ on the boundedness locus $\mathcal{M}_{f,c}$ and $h_{(f,c)}>0$ for any $t$ on the non-boundedness locus $\pran\setminus\mathcal{M}_{f,c}$ by definition and since we assume that the both sets are non-empty. We have thus the activity measure is a strictly positive  finite Radon measure. Indeed, If the activity measure $\mu_{(f,c)}$ is zero measure, the function $h_{(f,c)}$ should be constant since every harmonic function on $\affan$ is constant. 
\end{proof}
The relation between this set and the activity measure is the following:

\begin{thm}\label{slicemain}
\textit{ In the above notation, we assume $\mathcal{M}_{f,c}$ is non-empty and bounded. Then, the activity measure $\mu_{(f,c)}$ coincides with the equilibrium measure $\mu_{\mathcal{M}_{f,c}}$ of the set $\mathcal{M}_{f,c}$ with respect to $\infty$ times $\mu_{(f,c)}(V)$. Moreover, the support of the activity measure coincides with the boundary of $\mathcal{M}_{f,c}$.}
\end{thm}

\begin{rmk}\mbox{}\\
\vspace{-5mm}\begin{itemize}
\item For the definition of equilibrium measure, see Section 2. This is well-defined only when the set $\mathcal{M}_{f,c}$ is compact and of positive capacity. Since the positivity of the capacity is shown later, let us see the compactness of the set here. Since $\mathcal{M}_{f,c}$ is the zero locus of the dynamical Green function $h_{(f,c)}$ and it is continuous, it is closed. Also we assumed that the set is bounded. Hence it is also closed in the projective line $\pran$, which is compact and Hausdorff. 
\item As $\mu_{(f,c)}$ is positive and finite, we get a probability measure by taking a normalization $\mu_{(f,c)}/\mu_{(f,c)}(V)$. The above theorem states this probability measure coincides with the equilibrium measure of $\mathcal{M}_{f,c}$.
\end{itemize}
\end{rmk}

To prove the theorem, we need several lemmas. Before stating them, let us prepare some notations. Since the function $c_n$ is algebraic, it can be extended to a function $\overline{V}\to\pran$ where $\overline{V}=V\cup\{\infty\}$ as we define in the beginning of this section. We can also extend the function $h_{(f,c)}$ to a function defined on $\overline{V}$ by
\begin{align*}
	h_{(f,c)}(\infty)=\limsup_{x\to\infty, x\in V}h_{(f,c)}(x).
\end{align*}
We denote these extensions by $c_n$ and $h_{(f,c)}$, respectively, by abuse of notation. For the activity measure $\mu=\mu_{(f,c)}$, there exists a probability measure $\hat{\mu}$ defined on $\overline{V}$ such that
\begin{align*}
	\hat{\mu}|_{V}=\frac{\mu_{(f,c)}}{\mu_{(f,c)}(V)}.
\end{align*}

Indeed, since the extended dynamical Green function $h_{(f,c)}$ is of bounded differencial variation, we can consider the Laplacian of $h_{(f,c)}$, which gives us the formula
\begin{align*}
	\dd^c h_{(f,c)}=\delta_{\infty}-\hat{\mu} \text{ on } \overline{V}.
\end{align*}

\begin{lmm}\label{supp}
\textit{The support of $\hat{\mu}$ is contained in the boundary of the set $\mathcal{M}_{f,c}$.}
\end{lmm}
\begin{proof}
Note that $\Supp(\hat\mu)\subset \pran\setminus\mathcal{M}_{f,c}^{\circ}$ since $h_{(f,c)}\equiv 0$ on $\mathcal{M}_{f,c}$. Hence it is enough to show that for every point $t_0\in\pran\setminus\mathcal{M}_{f,c}$, there exists an open neighborhood $U\subset\mathcal{M}_{f,c}$ of $t$ such that $\hat{\mu}|_{U}$ is zero. Set $f_{t}(z)=\sum_{i=0}^{d} a_i(t)z^i$ with $a_d(t)\neq0$ and take a positive number $M$ such that $M>\max_{0\leq i\leq d-1}\{|a_i(t_0)|, 1\}$. Then, for every $z\in\pran$ with $|z|>M$, the canonical height of $f_{t_0}$ coincides with the naive one, i.e. 
\begin{equation}\label{bigheight}
h_{f_{t_0}}(z)=|z|.
\end{equation}
Indeed, for these $z$, we have $|z|^d>|a_i(t)z^i|$ for any $i=0,\ldots d-1$, so 
\begin{eqnarray*}
h_{f_{t_0}}(z)&=&\log\max\left(1, \left|\sum a_i(t)z^i\right|\right)\\
&=&\log\left|z^d\right|=d\log|z|
\end{eqnarray*}
from the ultrametric inequality.
Take an open neighborhood $U$ of $t_0$ so small that $M>\max_{0\leq i\leq d}\{|a_i(t)|\}$ for any $t\in U$. Then the equation (\ref{bigheight}) holds for all $t\in U$ if we replace $t_0$ by $t$. Since the set $\{f_t^n(c(t))\}$ is unbounded by the definition of $\mathcal{M}_{f,c}$, there exists a natural number $N_0$ such that $|f_{t_0}^{N_0}(c(t_0))|>M$. Taking a smaller open neighborhood of $t_0$ in $U$ if necessary, we can assume that for ever $t\in U$ we have $|f_t^{N_0}(c(t))|>M$. Hence together with (\ref{bigheight}) we have
\[
h_{(f,c)}(t)=\frac{1}{d^{N_0}}\log\max(1, |f_t^{N_0}(c(t))|).
\]
Since 
\[
\dd^ch_{(f,c)}=\frac{1}{d^{N_0}}\bigg(\deg(c_{N_0})\delta_{\infty}-c_{N_0}^*\delta_{\zeta_{0,1}}\bigg)
\]
on $U$ where $\deg c_{N_0}$ is a degree of $c_{N_0}: V\simeq \aff\to\aff$ as a polynomial, and we have $h_{(f,c)}(t_0)>0=h_(f,c)(t)$ for any $t\in c_n^{-1}\{\zeta_{0,1}\}$. Hence we can take smaller $U$ to get $\dd^c h_{(f,c)}=0$ on $U$.

\end{proof}

Therefore, we have $\hat{\mu}=\mu_{(f,c)}$ on $V$. Also, we have the following by the above lemma:

\begin{lmm}\textit{
The measure $\mu$ has continuous potentials, i.e., for any non-classical point $\zeta$, the potential function $u_{\zeta,\mu}$ is continuous.
}\end{lmm}
\begin{proof}

Note that the measure $\mu$ can be written as the Laplacian of $-h_{(f,c)}$ on $V$ and as the Laplacian of $u(x)=0$ on $\overline{V}\setminus\mathcal{M}_{f,c}$, both of which are continuous subharmonic function. The claim follows from \cite[Proposition8.65]{BR}

%\begin{prop}[\cite{BR}, Proposition 8.65]
%\textit{ Let $\mu$ be a positive measure on $\pran$ for which $-\mu$ is locally the Laplacian of a continuous subharmonic function. Then $\mu$ has continuous potentials.}
%\end{prop}

%Since $\mu$ is supported on $\mathcal{M}_{f,c}$ and on $V=\affan$, we have a continuous subharmonic function $h_{(f,c)}$ whose Laplacian gives $-\mu$. It remains to show that there is an open neighborhood $U$ of $\infty$ and a continuous subharmonic function $u$ such that $\dd^c u=-\mu$ on $U$, but we can take $U=\overline{V}\setminus\mathcal{M}_{f,c}$ and $u=0$.
\end{proof}

Now we prove Theorem \ref{slicemain}. Let us consider the Arakelov-Green function $g_{\hat{\mu}}(x,y)$ discussed in Section \ref{subsec:arakelovgreen}. This can be written explicitly by means of the dynamical Green function:
\begin{lmm} \label{lmm:potentialeq}\textit{
There exists a real constant $C$ such that
\begin{equation}\label{AGheight}
g_{\hat{\mu}}(x,y)=h_{(f,c)}(x)+h_{(f,c)}(y)-\log\delta(x,y)_{\infty}+C.
\end{equation}}
\end{lmm}
\begin{proof}
Take the Laplacian of each side with respect to $x$:
\begin{eqnarray*}
\dd^c_x g_{\hat{\mu}}(x,y)&=&\delta_{y}-\hat{\mu}, \mbox{ and}\\
\dd^c_x(h_{(f,c)}(x)+h_{(f,c)}(y)-\log\delta(x,y)_{\infty}&=&\delta_{\infty}-(\hat{\mu}-\delta_{\infty}-\delta_y)\\
&=&\delta_y-\hat{\mu}.
\end{eqnarray*}
Hence there exists a constant $C(y)$ depending on $y$ such that
\[
g_{\hat{\mu}}(x,y)=h_{(f,c)}(x)+h_{(f,c)}(y)-\log\delta(x,y)_{\infty}+C(y).
\]
Take the Laplacian with respect to $y$, we have that $C(y)$ is constant, which shows the lemma.
\end{proof}

The former part of Theorem \ref{slicemain} follows from the next lemma;
\begin{lmm}\textit{
    In the above settings, the two measures $\hat{\mu}$ and $\mu_{\mathcal{M}_{f,c}}$ coincide.
}\end{lmm}

\begin{proof}
We denote by $I_{\hat{\mu}}(\rho)$ the energy integral for a probability measure $\rho$:
\[
I_{\hat{\mu}}(\rho)=\iint_{\pran\times\pran} g_{\hat{\mu}}(x,y)d\rho(x)d\rho(y),
\]
which takes its minimum if and only if $\rho=\hat{\mu}$ by proposition \ref{EMP}. 

The definition of the equilibrium measure is the minimizer of another energy integral $I_{\infty}$; as in Section \ref{sec:capacity}, for a probability measure $\rho$ supported on $\mathcal{M}_{f,c}$,
\begin{align*}
    I_{\infty}(\rho)=\iint_{\mathcal{M}_{f,c}\times\mathcal{M}_{f,c}}-u_{\infty,y}(x)d\mu(x)d\mu(y)
\end{align*}
takes its minimum if and only if $\rho=\mu_{\mathcal{M}_{f,c}}$.

Now we compare the above two energy minimizing principle and Lemma \ref{lmm:potentialeq}.

Since the support of $\hat{\mu}$ is contained in $\partial\mathcal{M}_{f,c}$ by lemma \ref{supp}, we have $\hat{\mu}\in\mathcal{P}(\mathcal{M}_{f,c})$. Since $h(x)= 0$ on the set $\mathcal{M}_{f,c}$, we have by (\ref{AGheight})
\begin{eqnarray}
I_{\hat{\mu}}(\rho)&=&\iint g_{\hat{\mu}}(x,y)d\rho(x)\rho(y) \nonumber \\
&=&\iint (h(x)+h(y)+\log\delta(x,y)_{\infty}+C)d\rho(x)\rho(y) \nonumber \\
&=&\iint\log\delta(x,y)_{\infty}d\rho(x)\rho(y)+C=I_{\infty}(\rho)+C
\end{eqnarray}
for any $\rho\in\mathcal{P}(\mathcal{M}_{f,c})$.
since $\hat{\mu}\in\mathcal{P}(\mathcal{M}_{f,c})$, we see that $V_{\infty}(E)\leq I_{\hat{\mu}}(\rho)-C<\infty$ which shows the well-definedness of the measure $\mu_{\mathcal{M}_{f,c}}$. Moreover, $I_{\hat{\mu}}(\rho)$ takes its minimum if and only if $\rho=\hat{\mu}$, while $I_{\infty}(\rho)$ does if and only if $\rho=\mu_{\mathcal{M}_{f,c}}$. Therefore we get $\hat{\mu}=\mu_{\mathcal{M}_{f,c}}$.
\end{proof}

 Since $\hat{\mu}$ has continuous potentials, so is $\mu_{\mathcal{M}_{f,c}}$, which implies $\Supp\mu_{\mathcal{M}_{f,c}}=\partial\mathcal{M}_{f,c}$ by proposition \ref{contipot}. This concludes the proof of Theorem \ref{slicemain}.

\subsection{Boundedness condition for the boundedness locus $\mathcal{M}_{(f,c)}$} \label{bddcond}

In this subsection, we give a sufficient condition for the boundedness locus $\mathcal{M}_{(f,c)}$ to be bounded.

\begin{prop}\label{bddloci}
For $j=0,1,\ldots,d-2$, define $\mathcal{M}_{j}=\mathcal{M}_{f_{(\textbf{c},a)},c_j}$. Then,
\[
\overline{\M}_{j}\cap\mathbb{P}_{\infty}\subset H_{j}:=\overline{\bigg\{[\textbf{c}:a:0]\bigg|f_{(\textbf{c},a)}(c_{j}(\textbf{c},a))=0\bigg\}},
\]
where $\overline{E}$ is the topological closure of $E$ in $\mathbb{P}^{d-2, an}_{K}$ for a set $E$.
\end{prop}

Admitting the above proposition, we have a sufficient condition for $\mathcal{M}_j$ to be bounded; if $\overline{\mathcal{M}}_j$ does not intersect with $H_j$, the set $\mathcal{M}_j$ is bounded. Since $H_j$ is a hypersurface of $\mathbb{P}_{\infty}$, for a general choice of the line $\overline{V}$ of $\mathbb{P}$, the boundedness locus $\mathcal{M}_j$ is bounded.

To show the proposition, we need several lemmas;

\begin{lmm}\label{nonzerolmm}
$f_{(\textbf{c},a)}(c_{j})=0$ for any $j$ if and only if $(\textbf{c},a)=(0,\ldots,0,0)$.
\end{lmm}
\begin{proof}
First, as $f_{(\textbf{c},a)}(0)=0$, we have $a=0$. Note that for any $c_i$ with $f_{(\textbf{c},a)}(c_{i})=0$, the multiplicity $m(c_{i})$ is $1+($number of $j$ s.t. $c_{i}=c_{j})$ since the set of the critical points of $f_{(\textbf{c},a)}(z)$ is $\{c_0(=0),c_{1},\ldots, c_{d-2}\}$ and the multiplicity is the order of $f_{(\textbf{c},a)}(0)=0$ at $0$. Suppose there exists a non-zero critical point. Then, 
\[
d=\sum_{f_{(\textbf{c},a)}(z)=0}m(z)=\sum_{(*)} m(c_{j})\geq 2+(d-1),
\]
where (*) is the set of the distinct critical points. As this is contradiction, we have $\textbf{c}=(0,\ldots,0)$.
\end{proof}
Note that this is a purely algebraic result, so we have the same result for an arbitrary field.

For $(\textbf{c},a)\in K^{d-1}$, define following functions:
\begin{align*}
g_{(\textbf{c},a)}(z)&:=\lim_{n}\frac{1}{d^{n}}\log^{+}|f_{(\textbf{c},a)^{n}}(z)|,\\
g_{j}(\textbf{c},a)&:=g_{(\textbf{c},a)}(c_{j}), \text{ and}\\
G(\textbf{c},a)&:=\max\{g_{j}(\textbf{c},a)\},
\end{align*}
where $\log^+|z|=\log\max(1,|z|)$ for any $z\in K$.

\begin{lmm}\label{maxheight}
\[
G(\textbf{c},a)=\max\{|c_{j}|,|a|,1\}.
\]
\end{lmm}
\begin{proof}
Set $c_{j}=\pi^{\alpha_{j}}\cdot u_{j}$ for $j=1,\ldots d-2$ and $a=\pi^{\alpha_{0}}\cdot u_{0}$, where $\alpha_j$, $\alpha_0$ is a real number, $\pi$ is an element in $K$ such that $|\pi|<1$, and $|u_{j}|=1$.  As the claim is trivial when $\log^{+}\max\{|c_{j}|,|a|\}=1$, we assume that $\log^{+}\max\{|c_{j}|,|a|\}>1$. Let $I\subset\{0,\ldots,d-1\}$ be the set of all the indices that $|\pi^{\alpha_{j}}|$ attains the max in $G(\textbf{c},a)$. Then, $G(\textbf{c},a)$ is attained by some of $c_{j}$ with $j\in I$. Assume first $0\not\in I$. For any $j\in I$, we have
\begin{align*}
\left|f_{(\textbf{c},a)}(c_{j})\right|&=\left|\frac{1}{d}c_{j}^{d}+\sum_{i=1}^{d-2}\frac{(-1)^{d-i}\sigma_{i}(c)}{d-i}c_{j}^{d-i}+a^{d}\right|\\
&=\left|\pi^{d\alpha_{j}}\right|\cdot\left|\frac{1}{d}u_{j}^{d}+\sum_{i=1}^{d-2}\frac{(-1)^{d-i}\sigma_{i} (\{u_{k}\}_{k\in I},0,\ldots,0)}{d-i}u_{j}^{d-i}\right|,
\end{align*}
where $\sigma_{i} (\{u_{k}\}_{k\in I},0,\ldots,0)$ is the $i$-th fundamental symmetric polynomial of degree $d$ generated by $\{u_{k}\}_{k\in I}$ and $d-2-|I|$ times of 0's (considered as indeterminates). By Lemma \ref{nonzerolmm}, there exists $j\in I$ such that

\begin{equation}\label{nonzero}
\frac{1}{d}\tilde{u}_j^{d}+\sum_{i=1}^{d-2}\frac{(-1)^{d-i}\sigma_{i} (\{\tilde{u}_k\}_{k\in I},0,\ldots,0)}{d-i}\tilde{u}_j^{d-i}\neq0,
\end{equation}
where $\tilde{u}_k$ is the reduction of $u_k$ modulo the maximal ideal of $K$. For this $j$, we have $|f_{(\textbf{c},a)}(c_{j})|=|c_{j}|^{d}$.
Now, since it is easy to see that $|f_{(\textbf{c},a)}(z)|=|z|^{d}$ for any $z$ with $|z|>\max\{|c_{i}|,|a|\}$, We have 

\[
G(\textbf{c},a)=g_{j}(\textbf{c},a)=\frac{1}{d^{n}}\log\left|f_{(\textbf{c},a)}(c_{j})\right|=\log|c_{j}|=\max\{|c_{i}|,|a|\}.
\]
when $j\neq0$ is the non-zero solution of (\ref{nonzero}).
When $0\in I$ i.e.\ $|a|$ attains the maximum, we have 
\[
\left|f_{(\textbf{c},a)}(c_{0})\right|=\left|f_{(\textbf{c},a)}(0)\right|=|a|^{d}.
\]
In this case we have
\[
G(\textbf{c},a)=g_{0}(\textbf{c},a)=\frac{1}{d^{n}}\log\left|f_{(\textbf{c},a)}(0)\right|=\log|a|=\max\{|c_{i}|,|a|\}.
\]
\end{proof}

Now we prove Proposition \ref{bddloci}. First, it is easy to extend the result of Lemma \ref{maxheight} to the analytification of $\mathbb{A}^{d-1,\an}_{K}$ by continuity of each function. That is, the functions $g_{(\textbf{c},a)}(z)$, $g_{j}(\textbf{c},a)$, $G(\textbf{c},a)$ and $\log^{+}\max\{|c_{i}|,|a|\}$ can be extended to a continuous function on $\mathbb{A}^{d-1, \an}_{K}$ and we have
\begin{align}\label{bermaxheight}
G(\textbf{c},a)=\log^{+}\max\{|c_{j}|,|a|\}.
\end{align}
 Take any convergent sequence $\{\lambda_{n}\}$ in $\M$ and let $\lambda_{0}$ be the limit. 
 
 For the projective space $\mathbb{P}^{d-1}_{K}=\Proj K[T_{0},\ldots,T_{d-1}]$, we have a natural open covering by affine space $\{U_{i}\}_{i=0}^{d-1}$, where $U_{i}=\Spec K[T_{0}/T_{i},\ldots,T_{d-1}/T_{i}]$. Taking their analytification, we have an open covering $\{U_{i}^{\mathop{\mathrm{an
}}}\}$ of $\mathbb{P}^{d-1,\mathop{\mathrm{an}}}_{K}$. 

 Now, take one of $U_{i}^{\mathop{\mathrm{an}}}$ such that an infinitely many $\lambda_{n}$ belong to $U_{i}$ and satisfy $G(\lambda_{n})=g_{i}(\lambda_{n})$, and replace $\{\lambda_{n}\}$ by a subsequence with these properties. Without loss of generality, we may assume that $U:=U_{1}=\{[c_{1}:\cdots:c_{d-2}:a:k]|c_{1}\neq 0\}$, where $k=0$ on $\mathbb{P}_{\infty}$ and $k=1$ elsewhere and we can take an affine coordinate of $U$ as $(c_{2}/c_{1},\ldots,c_{d-2}/c_{1},a/c_{1},k/c_{1})$. Note that by assumption we have $g_{j}(\lambda)=0$ for any $n$ and $g_{1}(\lambda_{n})\to\infty$ as $n\to\infty$. By (\ref{bermaxheight}), we have
\begin{align*}
\log|c_{1}(\lambda_{n})|=F(\lambda_{n})=G(\lambda_{n})\geq\log\left|f_{\lambda_{n}}(c_{j}(\lambda_{n}))\right|,
\end{align*}
whence
\begin{equation}\label{maxandj}
\log|c_{1}(\lambda_{n})|\geq\log\left|f_{\lambda_{n}}(c_{j}(\lambda_{n}))\right|.
\end{equation}
Note that 
\[
f_{\lambda}(c_{j}(\lambda)/c_{1}(\lambda))=\frac{1}{c_{1}(\lambda)^{d}}f_{\lambda}(c_{j}(\lambda))
\]
since $f_{\lambda}(\lambda)$ is a homogeneous polynomial. It is enough to show $f_{\lambda}(c_{j}(\lambda_{n})/c_{1}(\lambda_{n}))\to 0$ as $n\to\infty$. From (\ref{maxandj}) we have
\[
(1-d)\log|c_{1}(\lambda_{n})|\geq\log\left|f_{\lambda_{n}}(c_{j}(\lambda_{n})/c_{1}(\lambda_{n}))\right|.
\]
Taking $n\to\infty$, we see that the left hand side tends to $-\infty$, which means $f_{\lambda}(c_{j}(\lambda_{n})/c_{1}(\lambda_{n}))\to 0$.

\section{Normality vs.\ Passivity and J-stability vs.\ Activity}\label{Normal}
Let $V$ be a smooth and strictly $K$-analytic curve, $f: V\times\pran\to\pran$ be an analytic family of rational functions of degree $d$, and $c:V\to\pran$ be a marked critical point. As explained in Section 1, there are two ways to forumlate the activity in non-archimedean dynamics. The other approach by means of normality, is our main interest in this section.

 As in Section \ref{subsection:normality}, the family of functions $\{c_n(t)=f_t^n(c(t))\}_n$ is said to be \textit{normal} around $t_0$ if for any subsequence $\{c_{n_j}\}$, there is an open neighborhood $U$ of $t_0$ and a continuous function $u$ on $U$ such that $c_{n_j}$ pointwise converges to $u$. In this section, we compare the normality locus in this sense with the complement of the activity locus discussed in Section \ref{escaperate}.

In general, we have the following property:
\begin{prop}
\textit{In the above notations, if $c$ is attracted by an attracting periodic point $w$ at $t_0$, $(f,c)$ is passive around $t_0$.}
\end{prop}
\begin{proof}
Replacing $f$ by $f^n$ where $n$ is the exact period of $w$ and taking an M\"obius transformation, we may assume that $w(t)=0$ is a fixed point. Take an open neighborhood $U$ of $t_0$ small enough that $0$ is attractive fixed point and $c$ is converging to $0$. Then clearly $\{c_n(t)\}$ is normal as it converges to a constant function $0$. Consider a homogeneous lift $F_t=(P_t,Q_t)$ of $f_t$ and a suitable lift $C_t: V\to\mathbb{A}^{2,\mathop{\mathrm{an}}}\setminus\{0\}$ of $c(t)$. By definition, the activity measure on $U$ is given by the Laplacian of the following function:
\[
h_{(f,c)}(t)= \lim_{n\to\infty}\frac{1}{d^n}\log\max\left(\left|P_t^{(n)}(C_t)\right|,\left|Q_t^{(n)}(C_t)\right|\right),
\]
where $(P_t^{(n)},Q_t^{(n)})$ is the $n$-th iteration of $F_t$. As $c(t)$ is converging to $0$,
\[
\lim_{n\to\infty}|f_t^{n}(c(t))|=\lim_{n\to\infty}\left|\frac{Q_t^{(n)}(C_t)}{P_t^{(n)}(C_t)}\right|=0.
\]
Therefore, for any $n$ large enough, $|P_t^{(n)}(C_t)|>|Q_t^{(n)}(C_t)|$. As we take $F_t$ and $C_t$ so that both $P_t^{(n)}(C_t)$ and $Q_t^{(n)}(C_t)$ are not zero at the same time and neither of them has a pole, 
\[
\log\max\left(\left|P_t^{(n)}(C_t)\right|,\left|Q_t^{(n)}(C_t)\right|\right)=\log\max\left|P^{(n)}_t(C_t)\right|
\]
is harmonic on $U$ for every $n$ large enough. Therefore, $h_{(f,c)}(t)$ is also harmonic on $U$ and $(f,c)$ is passive around $t_0$.
\end{proof}

For families of polynomials considered in Section \ref{poly}, we have the complete comparison as introduced in section 1:

\begin{prop}\label{main3}
\textit{Consider the same family as in the first paragraph of Section \ref{poly-1} and assume that its boundedness locus $\M_{(f,c)}$ is non-empty and bounded. then the passivity locus coincides with the normality locus.}
\end{prop}

\begin{proof}
First, suppose $t_0\in\Supp\mu_{(f,c)}$ and assume $\{c_n(t):=f_t^n(c(t))\}_n$ is normal around $t_0$. We know that $\partial\M_{(f,c)}=\Supp\mu_{(f,c)}$. As $\{c_n\}$ is normal, there exists a subsequence $\{c_{n_j}\}_j$ converging to a continuous function $\phi$ around $t_0$. Any open neighborhood $U$ of $t_0$ includes some point $s$ such that $c_n(s)\to\infty$ as $n\to\infty$. Therefore, by the continuity of $\phi$ we have $\phi(t_0)=\infty$, which is contradiction as $t_0\in\M_{(f,c)}$ i.e.\ the orbit of $t_0$ by $c_n$ must be bounded.

For the contrary, take any point $t_0\not\in\Supp\mu_{(f,c)}$. Since $\Supp\mu_{(f,c)}=\partial \M_{(f,c)}$, the point $t_0$ belongs to either $\affan\setminus\M_{(f,c)}$ or $\M_{(f,c)}^{\circ}$, where $\M_{(f,c)}^{\circ}$ is the (topological) interior of $\M_{(f,c)}$. In case $t_0\in\affan\setminus\M_{(f,c)}$, there exists an open neighborhood $U$ of $t_0$ so that every point $t$ of $U$ has unbounded orbit. Then, a family of functions $\{c_n(t)\}_n$ on $U$ is normal since $c_n(t)\to\infty\in\pran$ on $U$. In case $t_0\in\M_{(f,c)}^{\circ}$, there exists a basic open neighborhood $U$ of $t_0$ contained in $\M_{(f,c)}$. Since $\M_{(f,c)}$ is bounded, there is a real number $M$ such that $|c_n(t)|<M$ for every $t\in\M_{(f,c)}$. In particular, we have
\[
\bigcup_{n=1}^{\infty}c_n(U)\subset \pran\setminus\left\{t\in\pran:|t|>M\right\},
\]
where $\left\{t\in\pran:|t|>M\right\}$ is a generalized Berkovich open ball centered at the point $\infty$. Then Proposition \ref{FKTthm} asserts the normality of $\{c_n\}$ on $U$.

%Non-archimedean version of Montel's theorem in \cite{FKT}, therefore, asserts the normality; If there exists a point $t'\in U$ such that $\deg_{t'}(c_n(t))\to\infty$, one can show that any subsequence of $\{c_n\}_n$ admits a sub-subsequence converging to a constant by Proposition 3.1 of \cite{FKT}. If $\deg_{t'}(c_n)$ is bounded for all $t'\in U$, Theorem 2.1 in \cite{FKT} shows the normality around $t_0$.

\end{proof}

\begin{cor}\label{main4}
{\it Assume the residue characteristic of $K$ is greater than $d$. In the above setting, every parameter $t_0$ where $f$ has an unstably indifferent periodic point is in the activity locus of some critical point.}
\end{cor}

\begin{proof}
We use the following proposition:
\begin{prop}[\cite{BIJL}, Theorem 1.2]\label{BIJL}
\textit{If the residue characteristic of $K$ is greater than $d$, then, for any rational function $\phi\in K(T)$ and its attracting periodic point $z$, there exists some critical point $c$ such that $\phi^n(c)$ is strictly attracted by $z$.}
\end{prop}
As $f_{t_0}$ has an unstably indifferent periodic point, for any sufficiently small neighborhood $U$ of $t_0$, there is a marked periodic point $w(t):U\to \pran$ of $f_t$. Moreover, there is some tangent vector $\vec{\zeta}\in T_{t_0}(U)$ such that $w(t)$ is attracting for every $t$ represented by $\vec{\zeta}$. Taking $U$ smaller if necessary, there exists a critical point $c$ such that for every type I point $t$ represented by $\vec{\zeta}$, $c(t)$ is strictly attracted by $w(t)$. Meanwhile, There exists a tangent vector $\vec{\xi}\in T_{t_0}(U)$ such that $w(t)$ is repelling for all $t$ represented by $\vec{\xi}$ i.e.\ for such $t$, $c(t)$ is never attracted by $w(t)$. Since any neighborhood of $t_0$ contains both of a point representing $\vec{\zeta}$ and one representing $\vec{\xi}$, $\{c_n:=f_t^n(c(t))\}_n$ cannot be normal around $t_0$.
\end{proof}
\begin{rmk}
Even if the residue characteristic does not exceed $d$, the above corollary is still valid as long as $f_t$ is tame for every $t$. Actually, the exact condition that a critical point is attracted by an fixed point if
\[
0<\lambda<|\deg_{\zeta,\vec{v}}(f_t)|^d
\]
for all $\zeta\in\{t\}\times\pran$ and every tangent vector $\vec{v}$ at $\zeta$, where $\lambda$ is a multiplier of the fixed point. Here, we say a rational function $\phi$ is tame if for every $z$ and its tangent vector $\vec{v}$, $\deg_{\zeta, \vec{v}}(\phi)$ cannot be divided by the residue characteristic of $K$. When the family is tame, we have $|\deg_{\zeta,\vec{v}}(f_t)|^d=1$, which means as long as a periodic point is attracting, there exists some critical point attracted by it.

However, it is not true when a family has wild ramification; see the next section for a counter example.
\end{rmk}

\section{Example: quadratic polynomials} \label{quad}
Now let us consider a concrete example. Assume the residue characteristic of $K$ is $0$ or greater than $2$. Let $f_t(z)=z^2+t$, $c(t)=0$, and consider the Mandelbrot set $\mathcal{M}$. For $|t|\leq 1$, $f_t$ has a good reduction and the orbit of $0$ is bounded. For $|t|>1$, we have $|c_n(t)|=|t|^{2^{n-1}}$. Indeed, $|f_t(0)|=|t|$ and for $n>0$ we have
\[
|c_n(t)|=\max((|t|^{2^{n-2}})^2, |t|)=|t|^{2^{n-1}}
\]
since $|t|^2>|t|>1$. This implies the orbit is unbounded. Hence $\mathcal{M}=\mathcal{D}(0,1)$ and the boundary $\partial\mathcal{M}=\zeta_{0,1}$. Since the boundary consists of a single point, the activity measure should be a Dirac mass at the Gauss point $\zeta_{0,1}$. This is actually true because $c_n(t)=f_t^n(c(t))$ has good reduction for every $n$ as a polynomial. We have the following fact:

\begin{prop}[\cite{FRL1}, Th\'eor\`eme E]
\textit{ A rational map $\phi$ of degree at least $2$ has good reduction if and only if the canonial measure of $\phi$ is a Dirac mass $\delta_{\zeta_{0,1}}$ at the Gauss point.}
\end{prop}
Hence, the canonical measure must be a Dirac mass at the Gauss point $\delta_{\zeta_{0,1}}$, so
\[c_n^*\delta_{\zeta_{0,1}}=2^n\delta_{\zeta_{0,1}},
\]
from which we have $\mu=\hat{\mu}=\delta_{\zeta_{0,1}}$.

On the other hand, the Julia set of $f_t$ with $|t|>1$ is a Cantor set in the annulus $\{ z\in \pr(k)|$ $|z|=|t|\}$ while that of $f_t$ with $|t|\leq 1$ is a singleton $\{\zeta_{0,1}\}$ since $f_t$ has good reduction. As a consequence a drastic bifurcation occurs at the point $\zeta_{0,1}$. Also \cite{TS} shows that the bifurcation locus equals the set of parameters where there exist some unstably indefferent periodic point or some repelling periodic point with multiplicity greater than $1$. The point $\zeta_{0,1}$ corresponds to the one with unstably indefferent periodic point. Indeed, there is only one attracting periodic point $z_0$ of $f_t$: $z_0=1-\sqrt{1-4t}$. Then, this is unstably indefferent if $t=\zeta_{0,1}$ since the absolute value of this is the following:
\begin{eqnarray*}
\bigg|1-\sqrt{1-4t}\bigg|=\begin{cases}
|t| \mbox{ if } |t|\leq 1\mbox{, and}\\
|t|^{1/2}\mbox{ otherwise}.
\end{cases}
\end{eqnarray*}

\begin{rmk}
In case a family $f$ has wild ramification for some $t$, the activity of critical points cannot describe the J-stability. We can easily see this in the same example as above but the residue characteristic is $2$. The proof is straightforward:
\begin{prop}
\textit{ Assume that the characteristic of the residue field of $K$ is $2$. Set $f_t(z)=z^2+t$, $c(t)=0$ and $c_n(t)=f_t^n(c(t))=f_t^n(0)$. Denote by $\mathcal{M}$ the Mandelbrot set as above and $\mathcal{B}(f)$ be the set of points in the parameter space where $f_t$ has unstably indifferent cycle. Then,
\[
\partial\mathcal{M}=\zeta_{0,1}\mbox{, while}
\]\[
\mathcal{B}(f)=\zeta_{0,4}.
\]
In this case, the above two sets do not coincide.}
\end{prop}
The argument which shows the Mandelbrot set is $\mathcal{D}(0,1)$ is valid in this case, too. Only thing to consider here is about $\mathcal{B}(f)$.
As we change the coordinate as $z\mapsto z-\sqrt{t}$, we have the family $\{g_t(z)=z-2\sqrt{t}z\}$. This has good reduction if $|t|\leq4$. Hence, when $|t|\leq4$, the polynomial $f_t$ has potentially good reduction. In this case, the Julia set $J(f_t)$ moves continuously, so this should be J-stable even though this motion is out of the scope of the definition of J-stability in \cite{TS} since the stability is defined only if the Julia set has some classical point. For this case, the same drastic bifurcation described above occurs at $\zeta_{0,4}$, instead of $\zeta_{0,1}$, different from the boundary of the Mandelbrot set. The problem is from the constant $\varepsilon$ in Theorem \ref{BIJL}. In this case, we have
\[
\varepsilon=\max(1, |2|^2)=1/4<1.
\]
To associate attracting periodic point to the behavior of critical points, we need to assume $\sigma=1$.
\end{rmk}
\begin{rmk}
In \cite{S}, it has been discussed that the non-existence of unstably indifferent cycle is not enough for the J-stablity to hold; in fact, one has to add the non-existence of type 1 repelling point with multiplicity greater than $1$. However, our activity measure cannot detect over which parameter there exists a type 1 repelling point with multiplicity greater greater than 1, since this has nothing to do with the asymptotic behavior of the critical point. The bifurcation of this kind occurs, for instance, at every point of the open segment $(\zeta_{0,1},\infty)$ in the above example, on which the critical point $0$ escapes to $\infty$. In this case, the Julia set is always homeomorphic to a Cantor set, but there is no analytic motion of repelling-periodic points (for details, see \cite{S}). Our study via the activity measure, therefore, cannot detect the bifurcation phenomena of this type. However, we think it is still useful in studying bifurcations that changes the shape of Julia sets.

%In \cite{S}, it has been shown that the non-existence of unstably indifferent cycle is not enough for the J-stability to hold. In fact, the J-stability in non-archimedean dynamics is equivalent to the non-existence of the unstably indifferent cycle and of the type 1 repelling point with multiplicity greater then $1$.

%The equivalent condition of J-stability in non-archimedean dynamics from \cite{S} is not just the non-existence of unstably indefferent cycle, but also it requires the non-existence of the type 1 repelling point with multiplicity greater then $1$. However this activity measure cannot detect the point where there exists the type 1 repelling point with multiplicity greater than $1$ because it has nothing to do with the asymptotic behavior of the critical orbit. In the above example, this bifurcation occurs at every point on the open segment $(\zeta_{0,1},\infty)$, on which the critical point $0$ escapes to $\infty$. In this case, the Julia set is always homeomorphic to a Cantor set but we cannot take an analytic motion of repelling periodic points (for details, see \cite{S}).

%Hence the bifurcation we can study via the activity measure is limited.
\end{rmk}

\end{document}